\documentclass[11pt]{article}

\usepackage{hyperref}
\usepackage{ucs}
\usepackage[utf8]{inputenc}

\usepackage{threeparttable}

\usepackage{dutchcal}
\DeclareSymbolFontAlphabet{\mathcalorig}{symbols}

\usepackage{amsmath}
\usepackage{amssymb}
\usepackage{amsthm}
\usepackage[english]{babel}
\usepackage{fontenc}
\usepackage{graphicx, tikz}
\usepackage{tikz-cd}
\usetikzlibrary{arrows}
\usepackage{amscd}
\usepackage{microtype}
\usepackage{enumitem}

\usepackage{geometry}
 \usepackage{xcolor}
\hypersetup{
    colorlinks,
    linkcolor={red!50!black},
    citecolor={blue!50!black},
    urlcolor={blue!80!black}
}

\usepackage{listings}
\usepackage{color}

\definecolor{mygreen}{rgb}{0,0.6,0}
\definecolor{mygray}{rgb}{0.5,0.5,0.5}
\definecolor{mymauve}{rgb}{0.58,0,0.82}
\definecolor{backcolour}{rgb}{0.95,0.95,0.92}

\title{admcycles - a {S}age package for calculations in the tautological ring of the moduli space of stable curves}
\author{Vincent Delecroix, Johannes Schmitt, Jason van Zelm}

\date{\today}

 \newcommand{\todo}[1]{}
\newcommand{\todoOld}[1]{}
\newcommand{\todoAlt}[1]{}
\newcommand{\todoFin}[1]{}
\newcommand{\M}{\ensuremath{\overline{\mathcalorig{M}}}}
\newcommand{\DR}{\mathrm{DR}}

\newcommand{\comment}[1]{}

\newcommand{\detex}[1]{}  

\newcommand{\Sage}{SageMath}
\newcommand{\cocalc}{CoCalc}
\newcommand{\sagecell}{SageMathCell}
\newcommand{\admcycles}{\texttt{admcycles}}

\begin{document}

\maketitle

 \newtheoremstyle{test}
  {}
  {}
  {\it}
  {}
  {\bfseries}
  {.}
  { }
  {}
 
 \theoremstyle{test}
\newtheorem{Def}{Definition}[section]
\newtheorem{Exa}[Def]{Example}
\newtheorem{Rmk}[Def]{Remark}
\newtheorem{Exe}[Def]{Exercise}
\newtheorem{Theo}[Def]{Theorem}
\newtheorem{Lem}[Def]{Lemma}
\newtheorem{Cor}[Def]{Corollary}
\newtheorem{Pro}[Def]{Proposition}
\newtheorem*{Rmk*}{Remark}   
\newtheorem*{Exa*}{Example}   
\newtheorem*{Pro*}{Proposition} 
\newtheorem*{Def*}{Definition}
\newtheorem*{Cor*}{Corollary}
\newtheorem*{Lem*}{Lemma}
\newtheorem*{Theo*}{Theorem}

\begin{abstract}
 The tautological ring of the moduli space of stable curves has been studied extensively in the last decades. We present a \Sage{} implementation of many core features of this ring. This includes lists of generators and their products, intersection numbers and verification of tautological relations. Maps between tautological rings induced by functoriality, that is pushforwards and pullbacks under gluing and forgetful maps, are implemented. Furthermore, many interesting cycle classes, such as the double ramification cycles, strata of $k$-differentials and hyperelliptic or bielliptic cycles are available. In this paper we show how to apply the package, including concrete example computations.
\end{abstract}

\section{Introduction}
A crucial tool in the study of the singular cohomology of the moduli space $\M_{g,n}$ of stable curves is the tautological ring \[RH^*(\M_{g,n}) \subset H^*(\M_{g,n}) = H^*(\M_{g,n}, \mathbb{Q}).\]
It is a subring of the singular cohomology of $\M_{g,n}$ with an explicit, finite set of generators (indexed by decorated graphs $[\Gamma, \alpha]$) admitting combinatorial descriptions of operations like cup products and intersection numbers. 
For a detailed introduction to the tautological ring, see e.g. \cite{faber2000,Arbarello2011,calcmodcurves}.

Since computations with the generators $[\Gamma, \alpha]$ quickly become untractable by hand, it is natural to implement them in a computer program. With \admcycles{} we present such an implementation using the open source mathematical software \Sage{}~\cite{sage}. It is based on an earlier implementation by Aaron Pixton. It features intersection products and numbers between the classes $[\Gamma, \alpha]$ and verification of linear relations between these generators using the known generalized Faber-Zagier relations \cite{pixtonrels,pandhapixton,jandarels}. For the gluing and forgetful morphisms between (products of) the moduli spaces $\M_{g,n}$ it implements pullbacks and pushforwards of the generators $[\Gamma, \alpha]$ of the tautological ring.

Many geometric constructions of cohomology classes on $\M_{g,n}$ (such as the Chern classes $\lambda_d$ of the Hodge bundle $\mathbb{E}$ over $\M_{g,n}$) give classes contained in the tautological ring and can thus be written as linear combinations of classes $[\Gamma, \alpha]$. For many examples of such classes, the package \admcycles{} implements known formulas or algorithms to calculate them and thus allows further computations, such as intersections or comparisons to other cohomology classes. In particular, \admcycles{} contains
\begin{itemize}
    \item a formula for double ramification cycles $\DR_g(A)$ from \cite{Janda2016Double-ramifica} ,
    \item a conjectural formula for the strata $\overline{\mathcalorig{H}}_g^k(\textbf{m})$ of $k$-differentials from \cite{FP,SchmittDimension},
    \item (generalized) lambda classes, the Chern classes of derived pushforwards $R^\bullet \pi_* \mathcalorig{O}(D)$ of divisors $D$ on the universal curve $\pi: \mathcalorig{C}_{g,n} \to \M_{g,n}$, as discussed in \cite{PRvZ},
    \item admissible cover cycles\footnote{Computing these cycles was the original purpose of \texttt{admcycles}, hence the name of the package.}, such as the fundamental classes of loci of hyperelliptic or bielliptic curves with marked ramification points, as discussed in \cite{schmittvanzelm}.
\end{itemize}
Instead of discussing the details of the algorithms in \admcycles{}, this document serves as a user manual for the package, with an emphasis on concrete example computations. These computations are also available in an interactive online format on \cocalc{} (without need for registration):
\begin{center}
\href{https://share.cocalc.com/share/0a48957b67f375b9e3107216504ca0c4efb678fd/admcycles%20tutorial.ipynb?viewer=share}{Click here to open example computations on share.cocalc.com}.
\end{center}
One way to explore \admcycles{} is to go through these examples and refer back to the text below for additional explanations and background.

\subsection*{Applications of \admcycles{}}
By now the package \admcycles{} has been used in a variety of contexts. Its original purpose was computing new examples of admissible cover cycles in~\cite{schmittvanzelm}, e.g. computing the class of the hyperelliptic locus in $\overline{\mathcalorig{M}}_5$ and $\overline{\mathcalorig{M}}_6$ and the locus of bielliptic cycles in $\overline{\mathcalorig{M}}_4$. It was also used to verify results about Hodge integrals on bielliptic cycles in~\cite{hurwitzhodge} and on loci of cyclic triple covers of rational curves in~\cite{somerstep}.

The authors of~\cite{rossiburyak} used \admcycles{} to explore formulas for intersection numbers involving double ramification cycles and lambda classes. The implementation of generalized lambda classes led to the discovery of previously missing terms in the computations of~\cite{PRvZ} when doing comparisons with double ramification cycles. The package was also used in~\cite{2019arXiv191202267C} to verify computations of Masur-Veech volumes in terms of intersection numbers on $\overline{\mathcalorig{M}}_{g,n}$. Finally, it was used to check a new recursion for intersection numbers of $\psi$-classes presented in~\cite{2019arXiv190312526G}.

\subsection*{Other implementations}
Apart from \admcycles{} (and the code of Pixton on which it is based) there have been several other implementations of the tautological ring,  starting with the program~\cite{faberdivisors} by Faber for computing intersection numbers of divisors and Chern classes of the Hodge bundle. In~\cite{Yang2008} Yang presents a program computing intersection pairings of tautological classes on various open subsets of $\M_{g,n}$. The package~\cite{djohnson} by Johnson implements general intersections of the $[\Gamma, \alpha]$ and also verification of linear relations between these generators against the known generalized Faber-Zagier relations.

Based on the package~\admcycles{} there is work in progress by Matteo Costantini, Martin M\"oller and Jonathan Zachhuber on implementing tautological classes and intersection products on the smooth compactification of the strata of differentials presented in~\cite{BCGGM3}.





\subsection*{Acknowledgements}
We are indebted to Aaron Pixton for letting us use and modify a previous implementation of operations in the tautological ring by him as well as pointing out several issues in~\admcycles{} which have now been fixed. We thank Fr\'ed\'eric Chapoton, Samuel Leli\`{e}vre  and Jonathan Zachhuber for many valuable contributions.

We thank Harald Schilly and the team of \href{https://cocalc.com/}{CoCalc} as well as Andrey Novoseltsev and the team of \href{https://sagecell.sagemath.org/}{SageMathCell} for making \admcycles{} available on these platforms.

The first author was a guest of the Max-Planck Institut and then of the Hausdorff Institut for Mathematics during the development of the project.

The second author was
supported by the grant SNF-200020162928 and has received funding from the European Research Council (ERC)
under the European Union Horizon 2020 research and innovation programme
(grant agreement No 786580). During the last phase of the project, the second author profited from the SNF Early Postdoc.Mobility grant 184245 and also wants to thank the Max Planck Institute for Mathematics in Bonn for its hospitality. 

The third author was supported by the Einstein Foundation Berlin during the course of this work.

\subsection{Conventions} \label{Sect:conventions}
Let $\M_{g,n}$ be the moduli space of stable curves and $\pi: \M_{g,n+1} \to \M_{g,n}$ be the forgetful morphism of the marking $n+1$, which can be seen as the universal curve over $\M_{g,n}$. Let $\sigma_i : \M_{g,n} \to \M_{g,n+1}$ be the section of $\pi$ corresponding to the $i$th marked point ($i=1, \ldots, n$). For $\omega_\pi$ the relative dualizing line bundle of $\pi$ on the space $\M_{g,n+1}$ and $i=1, \ldots, n$ we define the $\psi$-class  
\[\psi_i = c_1( \sigma_i^* \omega_\pi) \in H^2(\M_{g,n}).\]
For $a=0,1,2, \ldots$ we define the (Arbarello-Cornalba) $\kappa$-class
\[\kappa_a = \pi_* \left( (\psi_{n+1})^{a+1} \right) \in H^{2a}(\M_{g,n}).\]
Finally, given a stable graph $\Gamma$ of genus $g$ with $n$ legs, let
\[\xi_\Gamma : \M_\Gamma = \prod_{v \in V(\Gamma)} \M_{g(v),n(v)} \to \M_{g,n}\]
be the gluing map associated to $\Gamma$. For a class $\alpha \in H^*(\M_\Gamma)$ given as a product of $\kappa$ and $\psi$-classes on the factors $\M_{g(v),n(v)}$, define
\[[\Gamma, \alpha] = (\xi_\Gamma)_* \alpha \in H^*(\M_{g,n}).\]
Such decorated boundary strata form a generating set (as a $\mathbb{Q}$-vector space) of the tautological ring $RH^*(\M_{g,n})$. 

\textbf{Note}: Throughout the paper and in the package \admcycles{}, we take the convention of \emph{not} dividing by the size $|\text{Aut}(\Gamma)|$ of the automorphism group of $\Gamma$. 
\tableofcontents

\section{Getting started}
The \admcycles{} package works on top of \Sage{} which is an open source software for mathematical computations. We describe how to install \Sage{} and \admcycles{} on a computer and how to use the available online services.

\subsection{\admcycles{} in the cloud}
The simplest way to play with \admcycles{} without installing anything beyond a web browser is to use one of \sagecell{} (\url{https://sagecell.sagemath.org/}, \cite{SageMathCell}) or the website \cocalc{} (\url{https://cocalc.com/}, \cite{CoCalc}). The former provides a basic interface to \Sage{} and \cocalc{}. The latter requires registration and allows to create worksheets that can easily be saved and shared. As mentioned before, it is possible to explore the computations presented below  \href{https://share.cocalc.com/share/0a48957b67f375b9e3107216504ca0c4efb678fd/admcycles%20tutorial.ipynb?viewer=share}{on share.cocalc.com} without the need to register.

\subsection{Obtaining \Sage{}}
\Sage{} is available on most operating systems. Depending on the situation one can find it in the list of softwares available from the package manager of the operating system. Alternatively, there are binaries available from the \Sage{} website~\url{https://www.sagemath.org/}. Lastly, one can compile it from the source code. More information on the installation process can be found at~\url{https://doc.sagemath.org/html/en/installation/}.

\subsection{Installation of the \admcycles{} package}
The package \texttt{admcycles} is available from the Python Package Index (PyPI) at
\begin{center}
    \url{https://pypi.org/project/admcycles/}
\end{center}
where detailed installation instructions are available for a range of systems. Note that the best performance (in particular for functions like \verb|DR_cycle|) is obtained using version 9.0 of \Sage{} or newer.

The package \admcycles{} is being developed on GitLab at
\begin{center}
\url{https://gitlab.com/jo314schmitt/admcycles}
\end{center}
where one can find the latest development version and a link to report bugs ("issues"). This is also the place to look at to suggest features or improvements.

\subsection{First step with \admcycles{}}
Once successfully installed, to use \admcycles{} one should start a \Sage{}-session and type
\begin{lstlisting}
sage: from admcycles import *
\end{lstlisting}
In the sample code, we reproduce the behavior of the \Sage{} console that provides the \texttt{sage:} prompt on each input line. When using the online \Sage-cell or a Jupyter worksheet, there is no need to write \texttt{sage:}. In all our examples, this \texttt{sage:} prompt allows to distinguish between the input (command) and the output (result).

\textbf{All other examples below assume that the line
\begin{center}
\texttt{from admycles import *} 
\end{center}
has been executed before.}

In addition to this manual, the package has an internal documentation with more informations concerning the various functions. To access additional information about some function or object \texttt{foo}, type \texttt{foo?} during the \Sage{} session.
\begin{lstlisting}
sage: psiclass?
Signature:      psiclass(i, g=None, n=None)
Docstring:     
   Returns the class psi_i on bar M_{g,n}.

   INPUT:

   i : integer
      The leg i associated to the psi class.

   g : integer
      Genus g of curves in bar M_{g,n}.

   n : integer
      Number of markings n of curves in bar M_{g,n}.

   EXAMPLES: ...
\end{lstlisting}

\section{Tautological classes} \label{Sect:Tautcl}
\subsection{Creating tautological classes} \label{Sect:entertaut}
There are different ways to enter tautological classes in the program and depending on the example, some are more convenient than others.
For the fundamental class, boundary divisors as well as $\psi$, $\kappa$ and $\lambda$-classes, there are predefined functions. 
\begin{itemize}
 \item \verb|fundclass(g,n)| returns the fundamental class of $\M_{g,n}$
 \item \verb|sepbdiv(h,A,g,n)| gives the pushforward $\xi_* [\M_\Gamma]$ of the boundary gluing map 
 \[\xi : \M_\Gamma = \M_{h,A \cup \{p\}} \times \M_{g-h,(\{1, \ldots, n\} \setminus A) \cup \{p'\}} \to \M_{g,n},\]
 where \verb|A| can be a list, set or tuple\footnote{Be careful that tuples of length $1$ must be entered as \texttt{(a,)} in Python, instead of \texttt{(a)}.} of numbers from $1$ to $n$. 
 \item \verb|irrbdiv(g,n)| gives the pushforward $(\xi')_* [\M_{g-1,n+2}]$ of the boundary gluing map 
 \[\xi' : \M_{g-1,n+2} \to \M_{g,n}\]
 identifying the last two markings to a node.  Note that, since $\xi'$ has degree $2$ onto its image, this gives \emph{twice} the fundamental class of the boundary divisor of irreducible nodal curves.
 \item \verb|psiclass(i,g,n)| gives the $\psi$-class $\psi_i$ of marking $i$ on $\M_{g,n}$
 \item \verb|kappaclass(a,g,n)| gives the (Arbarello-Cornalba) $\kappa$-class $\kappa_a$ on $\M_{g,n}$
 \item \verb|lambdaclass(d,g,n)| gives the class $\lambda_d$ on $\M_{g,n}$, defined as the $d$-th Chern class $\lambda_d=c_d(\mathbb{E})$ of the Hodge bundle $\mathbb{E}$, the vector bundle on $\M_{g,n}$ with fibre $H^0(C, \omega_C)$ over the point $(C,p_1, \ldots, p_n) \in \M_{g,n}$.
\end{itemize}

When working with a fixed space $\M_{g,n}$, it is furthermore possible to globally specify $g,n$ using the function \verb|reset_g_n(g,n)| to avoid giving them as an argument each time (see below).

These tautological classes can be combined in the usual way by operations \verb|+|, \verb|-|, \verb|*| and raising to a power \verb|^|.
\begin{lstlisting}
sage: t1=3*sepbdiv(1,(1,2),3,4)-psiclass(4,3,4)^2
sage: reset_g_n(2,1)
sage: t2=-1/3*irrbdiv()*lambdaclass(1)
\end{lstlisting}
To enter more complicated classes coming from decorated boundary strata, it is often convenient to first list all such decorated strata of a certain degree and then select the desired ones from the list. To give a list of all generators of $RH^{2r}(\M_{g,n})$ use the function \verb|list_tautgens(g,n,r)|. 
\begin{lstlisting}
sage: list_tautgens(2,0,2)
[0] : Graph :      [2] [[]] []
Polynomial : 1*(kappa_2^1 )_0 
[1] : Graph :      [2] [[]] []
Polynomial : 1*(kappa_1^2 )_0 
[2] : Graph :      [1, 1] [[2], [3]] [(2, 3)]
Polynomial : 1*(kappa_1^1 )_0 
[3] : Graph :      [1, 1] [[2], [3]] [(2, 3)]
Polynomial : 1*psi_2^1 
[4] : Graph :      [1] [[2, 3]] [(2, 3)]
Polynomial : 1*(kappa_1^1 )_0 
[5] : Graph :      [1] [[2, 3]] [(2, 3)]
Polynomial : 1*psi_2^1 
[6] : Graph :      [0, 1] [[3, 4, 5], [6]] [(3, 4), (5, 6)]
Polynomial : 1*
[7] : Graph :      [0] [[3, 4, 5, 6]] [(3, 4), (5, 6)]
Polynomial : 1*
\end{lstlisting}
The list itself is created by \verb|tautgens(g,n,r)|, from which one can then select the classes.
\begin{lstlisting}
sage: L=tautgens(2,0,2);
sage: t3=2*L[3]+L[4]
sage: t3

Graph :      [1, 1] [[2], [3]] [(2, 3)]
Polynomial : 2*psi_2^1 

Graph :      [1] [[2, 3]] [(2, 3)]
Polynomial : 1*(kappa_1^1 )_0
\end{lstlisting}
The output above should be interpreted as follows: each \verb|tautclass| consists of a sum of decorated boundary strata (represented by data type \verb|decstratum|), which consist of a graph (datatype \verb|StableGraph|) and a polynomial in $\kappa$ and $\psi$-classes (datatype \verb|kppoly|). 

\emph{Reminder:} For decorated stratum classes, we have the convention of \emph{not} dividing by the order of the automorphism group of the stable graph. Thus the elements of the list \verb|L| above are really just pushforwards of $\kappa$ and $\psi$ classes on products of moduli spaces $\M_{g(v),n(v)}$ under appropriate gluing maps.

Let us look at the example of generator \verb|L[3]| above. 
\begin{lstlisting}
[3] : Graph :      [1, 1] [[2], [3]] [(2, 3)]
Polynomial : 1*psi_2^1 
\end{lstlisting}
Its stable graph is represented by three lists.
\begin{enumerate}
 \item The first list \verb|[1, 1]| are the genera of the vertices, so there are two vertices, both of genus $1$. Note that vertices are numbered by $0,1,2, \ldots$, so in the above case, the vertices are numbers $0$ and $1$.
 \item the second list gives the legs (that is markings or half-edges) attached to the vertices, so vertex $0$ carries the half-edge $2$ and vertex $1$ the half-edge $3$,
 \item the third list gives the edges, that is half-edge pairs that are connected; in the above case, the two half-edges $2$ and $3$ form an edge, connecting the two vertices
\end{enumerate}
If we wanted to enter this \verb|StableGraph| manually, we could use its constructor as follows:
\begin{lstlisting}
sage: G = StableGraph([1,1],[[2],[3]],[(2,3)]); G
[1, 1] [[2], [3]] [(2, 3)]
\end{lstlisting}
        
The polynomial in $\kappa$ and $\psi$ is \verb|1*psi_2^1| in this case, so the half-edge $2$ on the second vertex carries a $\psi$-class. For the generator \verb|L[4]| the polynomial looks like \verb|1*(kappa_1^1 )_0|, meaning that vertex $0$ carries a class $\kappa_1^1=\kappa_1$. 

In Section \ref{Sect:bdrypushforward} we are going to see how to push forward tautological classes under any boundary gluing map. 

Finally, it is possible to manually enter tautological classes by using the constructors of the classes \verb|tautclass|, \verb|decstratum| and so on. We refer to the documentation of \verb|admcycles| for details on this. 

\subsection{Basic operations} \label{Sect:Basicop}
Apart from the usual arithmetic operations, we can take forgetful pushforwards and pullbacks of tautological classes and also compute the degree of tautological zero-cycles. In particular, we can compute intersection numbers. Below, for the forgetful map $\pi: \M_{1,3} \to \M_{1,2}$ forgetting the marking $3$ we verify the relations
\[\pi_* \psi_3^2 = \kappa_1, \pi^* \psi_2 = \psi_2 - D_{0,\{2,3\}}.\]
\begin{lstlisting}
sage: s1=psiclass(3,1,3)^2
sage: s1.forgetful_pushforward([3])

Graph :      [1] [[1, 2]] []
Polynomial : 1*(kappa_1^1 )_0 
sage: s2=psiclass(2,1,2)
sage: s2.forgetful_pullback([3])

Graph :      [1] [[1, 2, 3]] []
Polynomial : 1*psi_2^1 

Graph :      [1, 0] [[1, 4], [5, 3, 2]] [(4, 5)]
Polynomial : -1* 
\end{lstlisting}
Using the method \texttt{evaluate} of \texttt{tautclass}, we also compute the intersection number 
\[\langle \tau_0 \tau_1 \tau_2 \rangle_{1,3} = \int_{\M_{1,3}} \psi_1^0 \psi_2 \psi_3^2 = 1/12.\]
We check that it agrees with the prediction $\langle \tau_0 \tau_1 \tau_2 \rangle_{1,3} = \langle \tau_0 \tau_2 \rangle_{1,2} +\langle \tau_1 ^2 \rangle_{1,2}$ by the string equation.
\begin{lstlisting}
sage: s3=psiclass(2,1,3)*psiclass(3,1,3)^2
sage: s3.evaluate()
1/12
sage: s4=psiclass(2,1,2)^2+psiclass(1,1,2)*psiclass(2,1,2)
sage: s4.evaluate()
1/12
\end{lstlisting}
Note that in the current version of \texttt{admcycles}, the list of tautological generators $[\Gamma_i,\alpha_i]$ in a tautological class is \emph{not} automatically simplified by combining equivalent terms (since in general this requires testing graph isomorphisms between the $\Gamma_i$). When performing arithmetic operations with complicated tautological classes, this simplification can be manually triggered using the function \texttt{simplify}, as demonstrated in the following toy example:
\begin{lstlisting}
sage: psisum = psiclass(1,2,1) + 3 * psiclass(1,2,1); psisum
Graph :      [2] [[1]] []
Polynomial : 1*psi_1^1 

Graph :      [2] [[1]] []
Polynomial : 3*psi_1^1 
sage: psisimple = psisum.simplify(); psisimple
Graph :      [2] [[1]] []
Polynomial : 4*psi_1^1
\end{lstlisting}
In a future version of \texttt{admcycles} (after improving our algorithms for graph isomorphisms), we plan to automate this process.

\subsection{A basis of the tautological ring and tautological relations} \label{Sect:basis}
Using the function \verb|generating_indices(g,n,r)| one can compute the indices (for the list \texttt{tautgens(g,n,r)}) of a basis of $RH^{2r}(\M_{g,n})$, assuming that the generalized Faber-Zagier relations (see \cite{pixtonrels,pandhapixton,jandarels}) between the additive generators $[\Gamma, \alpha]$ give a complete set of relations between them. The function \verb|Tautvecttobasis| converts a vector with respect to the whole generating set into a vector in this basis. The function \verb|tautclass.toTautbasis(g,n,r)| converts a \verb|tautclass| into such a vector. \todo{expand FZconjectureholds}

Continuing the example from Section \ref{Sect:entertaut} we see:
\begin{lstlisting}
sage: generating_indices(2,0,2)
[0, 1]
sage: t3.toTautbasis(2,0,2)
(-48, 22)
\end{lstlisting}
This means that generators \verb|L[0]|, \verb|L[1]|  form a basis of $RH^4(\M_2)$ and the tautclass \verb|t3=2*L[3]+L[4]| is equivalent to \verb|-48*L[0]+22*L[1]|.

We can also directly verify tautological relations using the built-in function \verb|is_zero| of \verb|tautclass|. Below we verify the known relation $\kappa = \psi - \delta_0 \in R^1(\M_{1,n})$ for $n=4$. Here $\psi$ is the sum of all $\psi_i$ and $\delta_0$ is the sum of all separating boundary divisors, i.e. those having a genus $0$ component. For this, we list all stable graphs with one edge via \verb|list_strata(g,n,1)|. We exclude the graph \verb|gamma| with a self-loop by requiring that the number of vertices \verb|gamma.numvert()| is at least $2$. Then we can convert these graphs to tautclases by using \verb|to_tautclass|.
\begin{lstlisting}
sage: reset_g_n(1,4)
sage: bgraphs=[bd for bd in list_strata(g,n,1) if bd.numvert()>1]
sage: del0=sum([bd.to_tautclass() for bd in bgraphs])
sage: psisum=sum([psiclass(i) for i in range(1,n+1)])
sage: rel=kappaclass(1)-psisum+del0
sage: rel.is_zero()
True
\end{lstlisting}
It is also possible to express tautological classes in a basis of the tautological ring of suitable open subsets of $\M_{g,n}$, e.g. to verify that some relation holds on the locus of compact type curves. This works with the optional argument \texttt{moduli} 
of the functions \texttt{toTautbasis} and \verb|is_zero|. The following options are available for \texttt{moduli}:
\begin{itemize}
    \item \texttt{'st'} : all stable curves (default)
    \item \texttt{'tl'} : treelike curves (all cycles in the stable graph have length 1)
    \item \texttt{'ct'} : compact type (stable graph is a tree)
    \item \texttt{'rt'} : rational tails (there exists vertex of genus g)
    \item \texttt{'sm'} : smooth curves
\end{itemize}
For instance, we can verify the relation \[\lambda_1 = \frac{B_2}{2} \kappa_1 = \frac{1}{12} \kappa_1 \in H^2(\mathcalorig{M}_g)\] following from Mumford's computation \cite{mumfordtowards} in the case $g=3$:
\begin{lstlisting}
sage: kappaclass(1,3,0).toTautbasis(moduli='sm')
(1)
sage: lambdaclass(1,3,0).toTautbasis(moduli='sm')
(1/12)
\end{lstlisting}
In practice, much of the time in some computations is spent on calculating generalized Faber-Zagier relations between tautological cycles on $\M_{g,n}$. However, once computed, the relations can be saved to a file and be reloaded in a later session using the functions \verb|save_FZrels()| and \verb|load_FZrels()|. Careful: the function \verb|save_FZrels()| creates (and overwrites previous version of) a file \verb|new_geninddb.pkl| which depending on the previous computations can be quite large. 
\subsection{Pulling back tautological classes to the boundary} \label{Sect:tautpullback}
Below we pull back a generator of $RH^4(\M_{4})$ to the boundary divisor with genus partition $4=2+2$. This produces an element of type \verb|prodtautclass|, a tautological class on a product of moduli spaces, in this case $\M_{2,1} \times \M_{2,1}$. Two elements on the same product of spaces can be added and multiplied and further operations like pushforwards under (partial) gluing maps are supported. More details are given in the documentation of the class \texttt{prodtautclass}.

In the example below, we want to express the pullback to $\M_{2,1} \times \M_{2,1}$ in terms of a basis of $H^2(\M_{2,1} \times \M_{2,1})$ obtained from the preferred bases of the factors $H^*(\M_{2,1})$ given by \verb|generating_indices|. We can either represent the result as a list of matrices (giving the coefficients in the tensor product bases) or as a combined vector (using the option \verb|vecout=true|).

\begin{lstlisting}
sage: bdry=StableGraph([2,2],[[1],[2]],[(1,2)])
sage: generator=tautgens(4,0,2)[3]; generator
Graph :      [1, 3] [[2], [3]] [(2, 3)]
Polynomial : 1*psi_3^1
sage: pullback=bdry.boundary_pullback(generator); 
sage: pullback.totensorTautbasis(2)
[
                           [-3]
                           [ 1]
                  [0 0 0]  [-3]
                  [0 0 0]  [ 7]
[-3  1 -3  7  1], [0 0 0], [ 1]
]
sage: pullback.totensorTautbasis(2,vecout=true)
(-3, 1, -3, 7, 1, 0, 0, 0, 0, 0, 0, 0, 0, 0, -3, 1, -3, 7, 1)
\end{lstlisting}
%
%
\subsection{Pushing forward classes from the boundary} \label{Sect:bdrypushforward}
Recall that for a stable graph $\Gamma$ we have a gluing map
\[\xi_\Gamma : \M_{\Gamma} = \prod_{i=1}^m \M_{g(v_i),n(v_i)} \to \M_{g,n}\]
taking one stable curve for each of the vertices $v_1, \ldots, v_m$ of $\Gamma$ and gluing them together according to the edges of $\Gamma$. The pushforward under $\xi_\Gamma$ sends a product of tautological classes on the factors $\M_{g(v_i),n(v_i)}$ to a tautological class of $\M_{g,n}$. This operation is implemented by the function \verb|boundary_pushforward| of \texttt{StableGraph}.

That is, if \texttt{Gamma} is a \texttt{StableGraph} and \texttt{[c1, ..., cm]} is a list whose $i$th element \texttt{ci} is a \texttt{tautclass} on the $i$th factor $\M_{g(v_i),n(v_i)}$ of $\M_\Gamma$, then
\begin{center} \verb|Gamma.boundary_pushforward([c1, ..., cm])| \end{center} is the pushforward of the product of the \texttt{ci}. Here, the markings for the class \texttt{ci} are supposed to go from $1$ to $n(v_i)$, where the $j$th marking corresponds to leg number $j$ on the $i$th vertex of \texttt{Gamma}.

As an illustration, we verify that the package correctly computes the excess intersection formula proved in \cite{Graber2001} for the self-intersection of a boundary divisor in $\M_{3,3}$. 
\begin{lstlisting}
sage: B=StableGraph([2,1],[[4,1,2],[3,5]],[(4,5)])
sage: Bclass=B.boundary_pushforward() # class of undecorated boundary divisor
sage: si1=B.boundary_pushforward([fundclass(2,3),-psiclass(2,1,2)]); si1
Graph :      [2, 1] [[4, 1, 2], [3, 5]] [(4, 5)]
Polynomial : (-1)*psi_5^1 
sage: si2=B.boundary_pushforward([-psiclass(1,2,3),fundclass(1,2)]); si2
Graph :      [2, 1] [[4, 1, 2], [3, 5]] [(4, 5)]
Polynomial : (-1)*psi_4^1 
sage: (Bclass*Bclass-si1-si2).is_zero()
True
\end{lstlisting}
Note that eg. for the term \texttt{s2} we needed to hand the function the term \texttt{-psiclass(1,2,3)} in the first vertex, since in the graph \texttt{B} the half-edge \texttt{4} is leg number \texttt{1} in the list of legs at the first vertex (and we have $(g(v_1),n(v_1))=(2,3)$ for this vertex).

\section{Special cycle classes}
\subsection{Double ramification cycles} \label{Sect:DR}
A particularly interesting family of cycles on $\M_{g,n}$ is given by the double ramification cycles. Fixing $g,n$ they are indexed by nonnegative integers $k,d \geq 0$ and a tuple $A=(a_1,a_2, \ldots, a_n)$ of integers summing to $k(2g-2+n)$. 

The classical double ramification cycle (for $k=0,d=g$)
\[\DR_g(A) \in H^{2g}(\M_{g,n})\] 
has been defined as the pushforward of the virtual fundamental class of a space of maps to rubber $\mathbb{P}^1$ relative to $0, \infty$ with tangency conditions at $0, \infty$ specified by the vector $A$ (see \cite{liruan, Li2002A-degeneration-, Li2001Stable-morphism,Graber2005Relative-virtua}). In \cite{Janda2016Double-ramifica} it is shown  that this cycle is tautological and an explicit formula in terms of tautological generators is provided.

More precisely, for $g,n,k,d$ and $A$ with $A$ a partition of $k(2g-2+n)$, the paper constructs an explicit tautological class 
\[P_g^{d,r,k}(A) \in \mathbb{Q}[r] \otimes_{\mathbb{Q}} RH^{2d}(\M_{g,n})\]
with coefficients being polynomials in a formal variable $r$. We obtain a usual tautological class $P_g^{d,k}(A) \in RH^{2d}(\M_{g,n})$ by setting $r=0$ in these polynomial coefficients. Then it is shown (\cite[Theorem 1]{Janda2016Double-ramifica}) that in the special case $k=0,d=g$, this gives a formula for the double ramification cycle
\[\DR_g(A) = 2^{-g} P_g^{g,k}(A).\]
While this demonstrates that the cycle $P_g^{d,k}(A)$ is useful for $k=0, d=g$, it has many interesting properties for other values of $k,d$:
\begin{itemize}
    \item for $k$ arbitrary and $d=1$, the restriction of $2^{-1} P_g^{1,k}(A)$ to the compact-type locus $\mathcalorig{M}_{g,n}^\textup{ct}$ gives the pullback of the theta divisor on the universal Jacobian $\mathcalorig{J}$ over $\mathcalorig{M}_{g,n}^\textup{ct}$ under the extension of the Abel-Jacobi section
    \[\mathcalorig{M}_{g,n} \to \mathcalorig{J}, (C,p_1, \ldots, p_n) \mapsto (\omega_C^\textup{log})^{\otimes k}(- \sum_{i=1}^n a_i p_i),\]
    see \cite{Hain2013Normal-function, Grushevsky2012The-zero-sectio}.
    \item for $k$ arbitrary and $d=g$, various geometric definitions of a double ramification cycle have been put forward and an equality with $2^{-g} P_g^{g,k}(A)$ is conjectured (see \cite[Section 1.6]{holmesschmitt} for an overview). Recently, this equality has been established on $\mathcalorig{M}_{g,n}^\textup{ct}$ for $k \geq 1$ and one of the $a_i$ satisfying $a_i<k$ or $k \nmid a_i$, see \cite[Corollary 1.2]{holmesschmitt}.
    \item for $k$ arbitrary and $d>g$, the class $P_g^{d,k}(A)$ vanishes by \cite{claderjanda}.
\end{itemize}
In \texttt{admcycles}, the formula for $P_g^{d,k}(A)$ has been implemented. The function \verb|DR_cycle(g,A,d,k)| returns the cycle $2^{-d} P_g^{d,k}(A)$. The factor $2^{-d}$ was chosen such that \verb|DR_cycle(g,A)| indeed gives the cycle $\DR_g(A)$. With the option \texttt{rpoly=True}, it is even possible to compute the cycle $2^{-d} P_g^{d,r,k}(A)$ whose coefficients are polynomials in the variable $r$.

As an application, we can verify the result from \cite{HPS} that DR cycles 
satisfy the multiplicativity property 
\[\DR_g(A) \cdot \DR_g(B) = \DR_g(A) \cdot \DR_g(A+B) \in H^{4g}(\mathcalorig{M}_{g,n}^{tl})\]
on the locus $\mathcalorig{M}_{g,n}^{tl}$ of treelike curves but \emph{not} on the locus of all stable curves, in the example given in \cite[Section 8]{HPS}.
\begin{lstlisting}
sage: A=vector((2,4,-6)); B=vector((-3,-1,4))
sage: diff = DR_cycle(1,A)*DR_cycle(1,B)-DR_cycle(1,A)*DR_cycle(1,A+B)
sage: diff.is_zero(moduli='tl')
True
sage: diff.is_zero(moduli='st')
False
\end{lstlisting}
In fact, using that the cycle $\DR_g(A)$ is polynomial in the entries of the vector $A$ (i.e. a tautological class with polynomial coefficients), we can check multiplicativity for all vectors $A,B$ in the case $g=1, n=3$. To gain access to the polynomial-valued DR cycle, we define a polynomial ring and call \verb|DR_cycle| with a vector $A$ having as coefficients the generators of this ring:
\begin{lstlisting}
sage: R.<a1,a2,a3,b1,b2,b3> = PolynomialRing(QQ,6)
sage: A = vector((a1,a2,a3)); B = vector((b1,b2,b3))
sage: diff = DR_cycle(1,A)*DR_cycle(1,B)-DR_cycle(1,A)*DR_cycle(1,A+B)
sage: diff.is_zero(moduli='tl')
True
\end{lstlisting}
As a second application, we can verify the formula from \cite[Theorem 2.1]{rossiburyak} for intersection numbers of two DR cycles with $\lambda_g$ on $\M_{g,3}$ in the case $g=1$:
\begin{lstlisting}
sage: intersect = DR_cycle(1,A)*DR_cycle(1,B)*lambdaclass(1,1,3)
sage: f = intersect.evaluate(); factor(f)
(1/216) * (a2*b1 - a3*b1 - a1*b2 + a3*b2 + a1*b3 - a2*b3)^2
sage: g = f.subs({a3:-a1-a2,b3:-b1-b2}); factor(g)
(1/24) * (a2*b1 - a1*b2)^2
\end{lstlisting}

\subsection{Strata of \texorpdfstring{$k$}{k}-differentials}
Let $g,n,k \geq 0$ with $2g-2+n>0$ and let $\textbf{m} = (m_1, \ldots, m_n) \in \mathbb{Z}^n$ with $\sum_i m_i = k(2g-2)$. Consider the subset
\[\mathcalorig{H}_g^k(\textbf{m}) = \left\{(C,p_1, \ldots, p_n) \in \mathcalorig{M}_{g,n} : \omega_C^{\otimes k}\left(\sum_{i=1}^n m_i p_i \right) \cong \mathcalorig{O}_C \right\} \subset \mathcalorig{M}_{g,n}.\]
Denote by $\overline{\mathcalorig{H}}_g^k(\textbf{m})$ the closure of $\mathcalorig{H}_g^k(\textbf{m})$ inside $\M_{g,n}$. Since the above equality of line bundles is equivalent to the existence of a meromorphic $k$-differential $\eta$ on $C$ with zeros and poles exactly at the points $p_i$ with multiplicities $m_i$, the subsets $\overline{\mathcalorig{H}}_g^k(\textbf{m})$ are called \emph{strata of $k$-differentials}. 

These strata are of interest in algebraic geometry, the theory of flat surfaces and Teichm\"uller dynamics and have been studied intensely in the past. Elements appearing in the boundary have been classified in \cite{BCGGM1, BCGGM2} and a smooth, modular compactification has been constructed in \cite{BCGGM3}. The dimension of $\overline{\mathcalorig{H}}_g^k(\textbf{m})$ depends on $k, \textbf{m}$ as follows (see e.g. \cite{FP, SchmittDimension}).
\begin{figure}[ht]
\begin{center}
\begin{threeparttable}[b]
    \begin{tabular}{c|c|c|c}
    codim & $\textbf{m}=0$ & $\textbf{m}=k \cdot \textbf{m}'$ for $\textbf{m}' \in \mathbb{Z}^n_{\geq 0}$ & $\textbf{m} \neq k \cdot \textbf{m}'$ for $\textbf{m}' \in \mathbb{Z}^n_{\geq 0}$ \\ \hline
    $k=0$ & 0 & 0 & g\\
    $k=1$ & 0\tnote{1} & $g-1$ & $g$\\
    $k>1$ & 0\tnote{1} & $g-1$ and\tnote{2}~ $g$ & $g$
    \end{tabular}

\begin{tiny}    
\begin{tablenotes}
\item [1] This forces $g=1$.
\item [2] The set $\overline{\mathcalorig{H}}_g^1(\textbf{m}') \subset \overline{\mathcalorig{H}}_g^k(\textbf{m})$ is a union of components of codimension $g-1$ in $\M_{g,n}$, with all other components of $\overline{\mathcalorig{H}}_g^k(\textbf{m})$ having pure codimension $g$.
\end{tablenotes}
\end{tiny}
\end{threeparttable}
    \end{center}
\caption{Dimension theory of $\overline{\mathcalorig{H}}_g^k(\textbf{m})$}
\end{figure}

For $k \geq 1$, the papers \cite{FP, SchmittDimension} present conjectural relations between the fundamental classes $[\overline{\mathcalorig{H}}_g^k(\textbf{m})]$ and the formulas for the double ramification cycles proposed by Pixton (see Section \ref{Sect:DR}). As explained in the papers, these conjectures can be used to recursively determine all cycles 
\begin{itemize}
    \item $[\overline{\mathcalorig{H}}_g^k(\textbf{m})] \in RH^{2g}(\M_{g,n})$ for $k \geq 1$ and $\textbf{m} \neq k \textbf{m}'$ for some $\textbf{m}' \in \mathbb{Z}_{\geq 0}^n$,
    \item $[\overline{\mathcalorig{H}}_g^1(\textbf{m})] \in RH^{2g-2}(\M_{g,n})$ for $k = 1$ and $\textbf{m}  \in \mathbb{Z}_{\geq 0}^n$.
\end{itemize}
These recursive algorithms have been implemented in the function \texttt{Strataclass(g,k,m)}, where as above \texttt{m} is a tuple of $n$ integers summing to $k(2g-2)$.

As a small application, we can check that the stratum class $[\overline{\mathcalorig{H}}_2^1((3,-1))]$ vanishes (the stratum is empty since by the residue theorem there can be no meromorphic differential with a single, simple pole). Also, the stratum $\overline{\mathcalorig{H}}_2^1((2))$ exactly equals the class of the locus of genus $2$ curves with a marked Weierstrass point, which can be computed by the function \texttt{Hyperell} (see below for details).
\begin{lstlisting}
sage: L=Strataclass(2,1,(3,-1)); L.is_zero()
True
sage: L=Strataclass(2,1,(2,)); (L-Hyperell(2,1)).is_zero()
True
\end{lstlisting}

\subsection{Generalized lambda classes}
Let $\pi \colon \mathcalorig{C}_{g,n}\to \M_{g,n}$ be the universal curve and assume $n\geq 1$. Every divisor of $\mathcalorig{C}_{g,n}$, up to pullback of divisors on $\M_{g,n}$, takes the form 
\[
D= l\tilde{K} + \sum^n_{p=1}d_p\sigma_p + \sum_{\substack{h\leq g,\\ 1\in S\subset [n]}}a_{h,S} C_{h,S}
\]
for some integers $l$, $d_p$, $a_{h,S}$. Here $\tilde{K}=c_1(\omega_\pi)$ is the first Chern class of the relative dualizing sheaf,  $\sigma_p$ is the class of the $p$th section and 
\[
C_{h,S} = \xi_* \left[\M_{h,S\cup\{\bullet\}} \times \M_{g-h,[n]\backslash S \cup \{\star,x\}}\right]\in \mathrm{CH}^1(\M_{g,[n]\cup\{x\}})=\mathrm{CH}^1(\mathcalorig{C}_{g,n}).
\]
In \cite{PRvZ} a formula is given for the chern character $\operatorname{ch}(R^\bullet\pi_*\mathcalorig{O} (D))$. This chern character can be computed up to degree $\verb|dmax|$ using  $\verb|generalized_chern_hodge(l,d,a,dmax,g,n)|$. It takes as input an integer $\verb|l|$, a list $\verb|d=[d1,...,dn]|$ of the integers $d_i$ and a list of triples $\verb|a=[[h1,S1,ahS1],...,[hn,Sn,ahSn]]|$ where the $\verb|ahSi|$ are the integers $a_{h,S}$ above (given in any order). It is enough to just include the triples $\verb|[h,S,ahS]|$ for which $a_{h,S}$ is nonzero. 


Using \verb|generalized_lambda(i,l,d,a,g,n)| the chern class $c_i(-R^\bullet \pi_* \mathcal{O}(D))$ can be computed directly. In particular when $l=1$ and the $d_i$ and  $a_{h,S}$ are zero, this equals the normal $\lambda$ class
\begin{lstlisting}
sage: g=3;n=1
sage: l=1;d=[0];a=[]
sage: s=lambdaclass(2,g,n)
sage: t=generalized_lambda(2,l,d,a,g,n)
sage: (s-t).is_zero()
True
\end{lstlisting}

Let $d_1,...,d_n$ be integers such that $\sum_{i=1}^n d_i$ is divisible by $2g-2$. Let $l=\sum_{i=1}^n d_i /(2g-2)$ and let $a_{h,S}$ be integers such that 
\[
 D(\phi)=l\tilde{K} + \sum d_i \sigma_i + \sum a_{h,S}(\phi) C_{h,S}
\]
is $\phi$-stable on the locus of stable curves with one node (for definitions see \cite{KassPagani1} or \cite{PRvZ}). For the shifted\footnote{This shift is due to the fact that the literature on double ramification cycles uses the "log-convention", i.e. the entries of the input sum to $l(2g-2+n)$.} vector $A=(d_1+l, \ldots, d_n+l)$, the paper \cite{HKP} proves an equality
\begin{equation}\label{eq:DRphi}
\textup{DR}_g(A)|_{U(\phi)} = c_g(-R^\bullet \pi_* D(\phi))|_{U(\phi)}
\end{equation}
on the largest open locus $U(\phi)\subset \M_{g,n}$ where the Abel-Jacobi section $s_{l,d}(\phi)\colon  \M_{g,n}  \dashrightarrow \mathcalorig{J}_{g,n}(\phi)$
extends to a morphism. In particular $U(\phi)$ always includes $\mathcalorig{M}_{g,n}^{\text{ct}}$ and equals $\M_{g,n}$ if and only if $l$, $d$ is trivial or  $l(2g-2)=0$ and $d=[0,...,\pm 1,...,\mp 1, ...,0]$. See \cite[Section 4.3]{PRvZ} for more details.

The function \verb|DR_phi(g,d)| computes $c_g(-R^\bullet \pi_* D(\phi))$. We can verify equality (\ref{eq:DRphi}). 
\begin{lstlisting}
sage: g=2;d=[1,-1]
sage: (DR_cycle(g,d)-DR_phi(g,d)).is_zero()
True
\end{lstlisting}
We also see that equality does not always hold over all of $\M_{g,n}$ but it does hold over $\mathcalorig{M}_{g,n}^\textup{ct}$
\begin{lstlisting}
sage: g=2;d=[2,-2]
sage: (DR_cycle(g,d)-DR_phi(g,d)).toTautbasis()
(12, -4, 14, 7, -40, -10, -14, -12, 28, -4, 6, -1, 4, 0)
sage: (DR_cycle(g,d)-DR_phi(g,d)).toTautbasis(moduli='ct')
(0, 0, 0, 0, 0)
\end{lstlisting}

\subsection{Admissible cover cycles} \label{Sect:admccycles}
\subsubsection{Hyperelliptic and bielliptic cycles}
Before we go into details of how to specify general admissible cover cycles, let us mention the important cases of hyperelliptic and bielliptic cycles.

Given $g,n,m \geq 0$ with $n \leq 2g+2$ and  $2g-2+n+2m >0$, we have the locus $\overline{H}_{g,n,2m} \subset \M_{g,n+2m}$ of curves $(C,p_1, \ldots, p_n, q_1, q_1', \ldots, q_{m}, q_{m}')$ such that $C$ is hyperelliptic with $p_1, \ldots, p_n$ fixed points of the hyperelliptic involution and the pairs $q_i, q_i'$ being exchanged by this involution. An analogous definition gives the locus $\overline{B}_{g,n,2m} \subset \M_{g,n+2m}$ of bielliptic curves with $n \leq 2g-2$ fixed points and $m$ pairs of points forming orbits under the bielliptic involution.

Then the fundamental class of the (reduced) loci $\overline{H}_{g,n,2m}$ and $\overline{B}_{g,n,2m}$ can (in many cases) be computed by the functions \verb|Hyperell(g,n,m)| and \verb|Biell(g,n,m)| of our program.

As an example, we compute the class $[\overline{H}_{3}] \in RH^2(\M_3)$ and verify that we obtain the known result \[[\overline{H}_{3}]=9\lambda - \delta_0 - 3\delta_1,\]
where $\delta_0$ is the class of the divisor of irreducible nodal curves and $\delta_1$ is the divisor of curves with a separating node between a genus $1$ and a genus $2$ component.
\begin{lstlisting}
sage: H=Hyperell(3,0,0)
sage: H.toTautbasis()
(3/4, -9/4, -1/8)
sage: reset_g_n(3,0)
sage: H2=9*lambdaclass(1)-(1/2)*irrbdiv()-3*sepbdiv(1,())
sage: H2.toTautbasis()
(3/4, -9/4, -1/8)
\end{lstlisting}
Here we need to divide \verb|irrbdiv()| by two, the degree of the corresponding gluing map.

\subsubsection{Creating and identifying general admissible cover cycles}
In general, an admissible cover cycle is specified by a genus, a group as well as monodromy data. Currently, intersections are only implemented for cyclic groups. Below we will study bielliptic curves in genus $2$, which are double covers of elliptic curves branched over two points. As a first step we enter the ramification data.
\begin{lstlisting}
sage: G=PermutationGroup([(1,2)])
sage: list(G)
[(), (1,2)]
sage: H=HurData(G,[G[1],G[1]])
\end{lstlisting}

The function \verb|HurData| takes the group $G$ as the first argument and as the second a list of group elements $\alpha \in G$, each of which corresponds to the G-orbit of some marking $p \in C$. Here $\alpha$ is a generator of the stabilizer of $p$ under the group action $G \curvearrowright C$, which gives the monodromy around $p$. In other words, the natural action of the stabilizer $G_p = \langle \alpha \rangle$ on a tangent vector $v \in T_p C$ is given by \[\alpha . v = \exp(2 \pi i/\mathrm{ord}(h)) v.\]

Thus in the example above, we have two markings, both with stabilizer generated by \verb|G[1]=(1,2)| which acts by multiplication of $-1$ on the tangent space.

To identify the admissible cover cycle (inside the moduli space $\M_{g,n}$ with $n$ the total number of marked points from the ramification data) in terms of tautological classes, one can use the function \verb|Hidentify|. 

It pulls back the admissible cover cycle to all boundary divisors and (recursively) identifies the pullback itself in terms of tautological classes. It compares this pullback to the pullback of a basis of the tautological ring. Often this pullback map is injective in cohomology such that one can then write the admissible cover cycle in terms of the basis using linear algebra. Sometimes, it is necessary to additionally intersect with some monomials in $\kappa$ and $\psi$-classes.

To apply \verb|Hidentify| one gives the genus and the monodromy data as arguments. 
The standard output format is an instance of the class \verb|tautclass|. For those users familiar with Aaron Pixton's implementation of the tautological ring, there is the option \verb|vecout=true| which returns instead a vector with respect to the generating set of the tautological ring provided by this program. 
\begin{lstlisting}
sage: vbeta=Hidentify(2,H,vecout=true)
sage: vector(vbeta)
(517/4, -33, 11/4, 243/4, -125/4, 15/2, 41/4, 125, 99/4, -41, -1137/4, -285/4, 0, 0, 0, 0, 0, 0, -57/8, -3/8, 0, 0, 0, 0, 0, 0, 0, 0, 0, 0, 0, 0, 0, 0, 0, 0, 0, 0, 0, 0, 0, 0, 0, 0, 0, 0, 0, 0, 0, 0, 0, 0, 0, 0, 0, 0, 0, 0, 0, 0, 0, 0, 0, 0, 0, 0, 0, 0, 0, 0, 0, 0, 0, 0, 0, 0, 0, 0, 0, 0, 0, 0, 0, 0, 0, 0, 0, 0, 0, 0, 0, 0, 0, 0, 0, 0, 0, 0, 0, 0, 0, 0, 0, 0, 0, 0, 0, 0, 0, 0, 0, 0, 0, 0, 0, 0, 0, 0, 0, 0, 0, 0, 0, 0, 0, 0, 0, 0, 0, 0, 0, 0, 0, 0, 0, 0, 0, 0, 0, 0, 0, 0, 0, 0, 0, 0, 0, 0, 0, 0, 0, 0, 0, 0, 0, 0, 0, 0, 0, 0, 0)
\end{lstlisting}
The output above means specifically, that inside $\M_{2,2}$ the locus of bielliptic curves with the two points fixed by the involution being the marked points is given (in the generating set \verb|gens=all_strata(2,3,(1,2))| produced by Pixton's program) as $517/4 \cdot \verb|gens[0]| - 33 \cdot \verb|gens[1]| + \ldots$.

If we instead wanted to have bielliptic curves with two marked fixed-points of the involution and one pair of markings that are exchanged by the involution, we would need to use the monodromy data
\begin{lstlisting}
sage: H2=HurData(G,[G[1],G[1],G[0]])
\end{lstlisting}
in which case \verb|Hidentify(2,H2)| would live inside $RH^{8}(\M_{2,4})$.

If we only want to remember a subset of the markings, we can use the optional parameter \verb|marking| to give this subset. For instance, the command \verb|Hidentify(2,H,markings=[])| would give the pushforward of \verb|Hidentify(2,H)| in $\M_{2,2}$ to the space $\M_2$ under the forgetful morphism (see also Section \ref{Sect:B2}).


\subsubsection{Example: Specifying and Identifying \texorpdfstring{$[\overline B_2]$}{[B 2]} by hand} \label{Sect:B2}
The locus $\overline B_2 \subset \M_{2}$ of bielliptic curves is a divisor. A bielliptic genus $2$ curve is ramified over two points. In the following we want to use the methods of the previous section to identify its cycle class. 

Now when treating admissible cover cycles in general, our program a priori handles the cycle where all possible ramification points are marked. In this case, this is the cycle $[\overline B_{2,2,0}] \in RH^6(\M_{2,2})$ of bielliptic curves $C$ with the two ramification points $p_1, p_2$ marked. By specifying \verb|markings=[]| when calling \verb|Hidentify|, we tell it to remember none of the markings, in other words to push forward under the map $\pi: \M_{2,2} \to \M_2$ forgetting the markings.
\begin{lstlisting}
sage: G=PermutationGroup([(1,2)])
sage: H=HurData(G,[G[1],G[1]])
sage: Biell=Hidentify(2,H,markings=[])
sage: Biell.toTautbasis(2,0,1)
(30, -9)
\end{lstlisting}
We want to compare the result with the known formulas for $[\overline B_2]$. For $\delta_0$ the class of the irreducible boundary of $\M_2$ and $\delta_1$ the class of the boundary divisor with genus-splitting $(1,1)$, it is known that $[\overline B_2]=\frac{3}{2} \delta_0 + 6 \delta_1$ (see \cite{faberbiell}). If we want to enter this combination of $\delta_0$ and $\delta_1$, we have to be careful about conventions, though: the corresponding gluing maps $\xi: \M_{1,2} \to \M_{2}$ and $\xi':\M_{1,1} \times \M_{1,1} \to \M_2$ both have degree $2$. This corresponds to the fact that the associated stable graphs both have an automorphism group of order $2$. Hence we have to divide by a factor of two and obtain:
\begin{lstlisting}
sage: reset_g_n(2,0)
sage: Biell2=3/4*irrbdiv()+ 3*sepbdiv(1,())
sage: Biell2.toTautbasis(2,0,1)
(15/2, -9/4)
\end{lstlisting}
We see that up to a factor of $4$ the two vectors \verb|(30, -9)| and \verb|(15/2, -9/4)| agree. Where does this factor come from? 

For this recall that the cycle \verb|Biell| above is equal to $\pi_* [\overline B_{2,2,0}]$. Since for the generic bielliptic curve $C$ there are two choices of orderings for marking $p_1, p_2$, this explains a factor of $2$. On the other hand, the hyperelliptic involution $\sigma: C \to C$ on $C$ exchanges $p_1$ and $p_2$. Thus $\sigma \in \operatorname{Aut}(C)$, but $\sigma \notin \operatorname{Aut}(C,p_1,p_2)$. This missing automorphism factor explains another factor of $2$ in the pushforward under $\pi$, so in fact $[\overline B_2] = \frac{1}{4} \pi_* [\overline B_{2,2,0}]$. 

Note that since the cycles of bielliptic loci are implemented via the function \verb|Biell|, we could have taken a shortcut above.
\begin{lstlisting}
sage: B=Biell(2,0,0); B.toTautbasis()
(15/2, -9/4)
\end{lstlisting}

As an application, we can check the Hurwitz-Hodge integral
\[\int_{[\overline B_{2,2,0}]} \lambda_2 \lambda_0 = \int_{\pi_* [\overline B_{2,2,0}]} \lambda_2 = \frac{1}{48} \]
predicted by \cite{hurwitzhodge}.
\begin{lstlisting}
sage: (Biell*lambdaclass(2,2,0)).evaluate()
1/48
\end{lstlisting}
The corresponding integrals for $g=3,4$ have also been verified like this, but the amount of time and memory needed grows drastically.

We can also check the Hurwitz-Hodge integral
\[\int_{[\overline{\mathcalorig{H}}_{2,\mathbb{Z}/3\mathbb{Z},((1,2,3)^2,(1,3,2)^2)}]} \lambda_1  = \frac{2}{9} \]
of $\lambda_1$ against the locus of genus $2$ curves admitting a cyclic triple cover of a genus $0$ curve with two points of ramification $(1,2,3) \in \mathbb{Z}/3\mathbb{Z}$ and two points of ramification $(1,3,2) \in \mathbb{Z}/3\mathbb{Z}$,
computed in \cite[Section 5]{somerstep}.
\begin{lstlisting}
sage: G = PermutationGroup([(1,2,3)]); list(G)
[(), (1,2,3), (1,3,2)]
sage: H = HurData(G,[G[1],G[1],G[2],G[2]]) #n=2, m=2
sage: t = Hidentify(2,H,markings=[])
sage: (t*lambdaclass(1,2,0)).evaluate()
2/9 
\end{lstlisting}
Note that while originally the cycle $[\overline{\mathcalorig{H}}_{2,\mathbb{Z}/3\mathbb{Z},((1,2,3)^2,(1,3,2)^2)}]$ lives in $\overline{\mathcalorig{M}}_{2,4}$, since we intersect with $\lambda_1$ which is a pullback from $\overline{\mathcalorig{M}}_2$ we can specify \texttt{markings=[]} above to compute the pushforward \texttt{t} of this cycle to $\overline{\mathcalorig{M}}_2$ before intersecting. This significantly reduces the necessary computation time.
\bibliographystyle{alpha}
\bibliography{manual}

\vspace{+16 pt}
\noindent Vincent Delecroix \\
\noindent Laboratoire Bordelais de Recherche en Informatique \\
\noindent Universit\'e de Bordeaux \\
\noindent \href{mailto:vincent.delecroix@u-bordeaux.fr}{vincent.delecroix@u-bordeaux.fr}

\vspace{+16 pt}
\noindent Johannes Schmitt \\
\noindent Mathematisches Institut \\
\noindent Universit\"at Bonn \\
\noindent \href{mailto:schmitt@math.uni-bonn.de}{schmitt@math.uni-bonn.de}

\vspace{+16 pt}
\noindent Jason van Zelm \\
\noindent Institut f\"{u}r Mathematik \\
\noindent Humboldt-Universit\"{a}t zu Berlin \\
\noindent \href{mailto:jasonvanzelm@outlook.com}{jasonvanzelm@outlook.com}

\end{document}